\documentclass[12pt,reqno]{amsart}


\newtheorem{proposition}{Proposition}[section]

\newtheorem{theorem}{Theorem}
\newtheorem{corollary}{Corollary}
\newtheorem{lemma}{Lemma}
\newenvironment{remark}{\noindent {\bf Remark:}}{}
\renewenvironment{proof}{\noindent {\bf Proof.}}{ \hfill\qed\\ }

\def\qed{\vrule height 5pt width 5 pt depth 0pt}
\def\glu{\!\!\!}

\def\e{\varepsilon}

\def\Y{\mathcal{Y}}
\def\Z{\mathcal{Z}}
\def\RR{{\mathbb R}}
\def\1{{{\mathit 1} \!\!\>\!\! I} }
\renewcommand{\phi}{\varphi}
\renewcommand{\liminf}{\mathop{{\underline {\hbox{{\rm lim}}}}}}
\renewcommand{\limsup}{\mathop{{\overline {\hbox{{\rm lim}}}}}}
\renewcommand{\eqref}[1]{{(\ref{#1})}}
\newcommand{\eqdef}{:=}

\begin{document}

\title{Recurrence, dimensions and Lyapunov exponents}
\author{B.~Saussol \and S.~Troubetzkoy \and S.~Vaienti}
\address{LAMFA / CNRS fre 2270, Universit\'e de Picardie Jules Verne, 
33, rue St Leu, F-80039 Amiens cedex 1, France}
\email{benoit.saussol@mathinfo.u-picardie.fr}
\urladdr{http://www.mathinfo.u-picardie.fr/saussol/}
\address{Centre de physique th\'eorique, Institut de math\'ematiques de Luminy et Federation de Recherches des Unites de Mathematique de Marseille, CNRS Luminy, Case 907, F-13288 Marseille Cedex 9, France.}
\email{serge@cpt.univ-mrs.fr}
\email{troubetz@iml.univ-mrs.fr}
\urladdr{http://iml.univ-mrs.fr/{\lower.7ex\hbox{\~{}}}troubetz/}
\address{Phymat, Universit\'e de Toulon et du Var, Centre de physique th\'eorique et Federation de Recherches des Unites de
Mathematique de Marseille, CNRS Luminy, Case 907, F-13288 Marseille Cedex 9, France}
\email{vaienti@cpt.univ-mrs.fr}
\urladdr{http://www.cpt.univ-mrs.fr/}
\begin{abstract}
We show that the Poincar\'e return time of a typical cylinder is at least
its length.  For one dimensional maps we express the Lyapunov exponent and 
dimension via return times.
\end{abstract}
\maketitle

\bibliographystyle{plain}

\section{Introduction} 

The statistical description of dynamical system has recently been
enriched by the study of recurrence and return times. 
For example the Poisson distribution
for the return and entrance time into a given set has been investigated 
(see \cite{hsv} and references therein), and the 
dimension like characterstics of invariant sets have been investigated
by means of recurrence \cite{a},\cite{psv}. 
In these two contexts a local quantity plays a fundamental role: the first
return of a set into itself, sometimes also called the Poincar\'e recurrence of the set. If $A \subset X$ is a measurable set of
a measurable (probability) dynamical systems $\{X,\beta, \mu, T\}$, the first return of the set $A$ is simply defined as:
$$
\tau(A) = \min\{n>0: \, T^n A\cap A\neq \emptyset\}.
$$
We will suppose in the following that $A$ is either the cylinder of order 
$n$ around the point $x\in X$ with respect to a measurable
partition $\mathcal{U}$,
i.e.~$A = U_n^x \in
\bigvee_{i=0}^{n-1}T^{-i}\mathcal{U}$, $x\in U_n^x$, or $A$ is a ball of
radius $r$ around $x$, $A=B_r(x)$. It has been shown in \cite{hsv} that the 
limits:
\begin{equation}\label{s1}
\liminf_{n\rightarrow\infty} \frac{\tau({U_n^x})}{n}
\quad \quad \hbox{ and } \quad \quad
\limsup_{n\rightarrow\infty} \frac{\tau({U_n^x})}{n}
\end{equation}
define $\mu$--almost everywhere invariant (subadditive) functions that 
control the asymptotic short returns into the sets
$U_n^x$. It is therefore important to have information on the value(s) of 
\eqref{s1}. The first result of this article is to show that
for a measurable dynamical system with
positive metric entropy, the $\liminf$ in \eqref{s1} is $\mu$--a.e. bigger or 
equal to~$1$. This estimate is a key bound need in proving the
exponential and Poisson statistics for return times, as pointed out in 
\cite{hsv} and \cite{hv}. 

For systems with zero metric entropy, examples are known where
$\liminf$ is positive and strictly smaller than one, for instance in
the case of Fibonacci rotations \cite{chv}
and where the $\liminf$ is identically equal to zero \cite{acs}. 
Our proof of the lower
bound (Theorem \ref{thm:tau/n}) is 
surprisingly easy. It uses the concept of Kolmogorov
complexity, which we recall briefly in the Appendix. 

Besides the interest in the asymptotic distribution of return and
entrance times, the limits in \eqref{s1} can be exploited in another 
direction, which, in \cite{lpv}, we proposed to call the
``thermodynamics of return times''.  We prove (Theorem 
\ref{thm:lower}) that for a large class of ergodic one-dimensional maps, 
the Lyapunov exponent can be estimated from the behavior of the first 
return times of a ball as the radius vanishes.

We turn to computing the Hausdorff dimension of the measure for the same 
class of one-dimensional maps. The first return of a set will be now replaced by 
another quantity which will denote with $\tau_{r}(x)$ and that is the first return 
of the point $x$ into its neighborhood $B_r(x)$. 
\begin{equation}\label{ss1}
\lim_{r\rightarrow  0^+}\frac{\log \tau_r(x)}{-\log r}
\end{equation}
In our final result (Theorem \ref{thm:dim}) we 
show that the limit in \eqref{ss1} exists $\mu$-a.e. and is equal to
the local dimension:
\begin{equation}\label{ss2}
d_\mu(x) := \lim_{r\rightarrow  0^+}\frac{\log \mu(B_r(x))}{\log r}
\end{equation}
For the class of one-dimensional maps considered 
is known that the local dimension is almost everywhere constant and
equal to the Hausdorff dimension of the measure 
and is also equal to the ratio between the metric entropy and the 
$\mu$--Lyapunov exponent \cite{h}.
A similar result in a multidimensional setting has been proved by Barreira 
and Saussol \cite{bs,bs2}, namely for the basic sets of Axiom-A
diffeomorphisms. Our proof is relatively simple and uses sharp comparison 
between balls and cylinders provided by Hofbauer
\cite{h}; nevertheless it is quite general since it covers maps with 
critical and parabolic points.

In conclusion these results are one of the first steps in establishing what 
we have already 
called thermodynamics of return times. In this
context a major role is played by the first returns of sets and points which often play the role of the  measure of
balls and cylinders according to the old suggestion given by Kac's theorem. This approach could be even more advantageous in
numerical and experimental investigations of dynamical systems as we showed in \cite{lpv}; further theoretical developments
in this direction appeared quite recently in \cite{hv} and \cite{acs}.

We will systematically use the same letter and
underscore (overscore) the names of pairs of functions defined via a
$\liminf$ ($\limsup$).
Because of this convention we will only write the one of the definitions 
of such pairs of functions.

\section{No small returns}
In this section we provide sharp lower bounds for the first return time, 
for cylinders of measurable partitions. 
These estimates are essential to compute the speed
of convergence to the exponential law of the first return time.
To prove this theorem we will use White's sharpening \cite{w,w1} of
a remarkable theorem by Brudno \cite{b} which links Kolmogorov complexity 
to entropy.
We state this theorem and give a quick introduction to
Kolmogorov complexity in the appendix.

\begin{theorem}\label{thm:tau/n}
Let $(T,X,\mu)$ be an ergodic measure preserving dynamical system. 
If $\zeta$ is a finite or
countable measurable partition with entropy $h_\mu(T|\zeta)$  strictly 
positive and $\zeta_n^x$ is the cylinder of length $n$ containing $x$,
then the lower rate
of return for cylinders is almost surely bigger or equal to 1, i.e.
$$
\underline{R}_\zeta(x)\eqdef\liminf_{n\to\infty} 
\frac 1n \tau(\zeta_n^x) \geq 1 
$$
for $\mu$-a.e. $x\in X$.
\end{theorem}

\begin{remark}\label{rem:1}
It is an easy exercise to show that if additionally $(T,X)$ satisfies
the specification property \cite{kh} then
$$\overline{R}_\zeta(x) \leq 1 
$$
for $\mu$-a.e. $x\in X$ and thus
$$
R_\zeta(x)\eqdef\lim_{n\to\infty} 
\frac 1n \tau(\zeta_n^x) = 1 
$$
for $\mu$-a.e. $x\in X$. Afraimovich {\it et.~al.} have shown that there
are specific examples of zero entropy maps for
which the conclusion of Theorem~\ref{thm:tau/n} is not true \cite{acs}.
\end{remark}

\begin{proof}
It is sufficient to prove the theorem for finite partitions, 
the case of countable $\zeta$ will follow easily.
More precisely, if $\zeta=\{B_1,B_2,\ldots\}$ is a countable partition, 
then for some $m<\infty$ the 
partition $\hat\zeta = \{B_1,B_2,\ldots,B_m, \cup_{l>m}B_l\}$ will have positive entropy. 
In addition, $\zeta$ is finer than $\hat\zeta$, hence $\tau(\zeta_n^x) \geq \tau(\hat\zeta_n^x)$.

Thus we assume that $\zeta$ is finite.  We claim that a cylinder 
$\zeta_n^x$ is completely determined by its first  $\tau := \tau(\zeta_n^x)$
symbols. To see this suppose that $\zeta_n^x = [x_0, \dots, x_{n-1}]$ and
that $y \in \zeta_n^x$ satisfies $T^{\tau}y \in \zeta_n^x$.  Let $j$ be the 
integer defined by $j \tau < n \le (j + 1) \tau.$  Since $T^{\tau}y \in 
\zeta_n^x$ we have $x_{i \tau }, \dots , x_{(i + 1)\tau - 1} = x_0, \dots ,
x_{\tau - 1}$ for $i = 0,1,\dots, j,$ proving the claim.

We will use the notion of Kolmogorov complexity to prove the theorem. All
the notations used here are defined in the appendix, more details can be
found in the references \cite{b,w}.  Let $N = \# \zeta$ and let
$K(w)$ be the Kolmogorov complexity of a finite word $w$ with  entries from
the alphabet $\{0,1,\dots,N-1\}$.
The partition $\zeta$ gives rise to the symbolic space $\Sigma_{\zeta}$ 
and a map
$\phi:X \to \Sigma_{\zeta}$ which is a semiconjugacy $\sigma \phi = \phi T.$
Let $\zeta_n^x$ be the word consisting of the first $n$--symbols of 
the sequence $\phi(x)$.  We use the notation $K_\zeta(\zeta_n^x)$ for the
complexity of $\zeta_n^x$ and define 
$K_{\zeta}(x,T) = \liminf_{n \to0} \frac{1}{n}K(\zeta_n^x)$.

Since $\zeta_n^x$ is determined by its first $\tau$ symbols 
the complexity of the sequence $\zeta_n^x$ is bounded by the complexity of defining
the first $\tau$ symbols plus the complexity of repeating these symbols up to
the size $n$ of the cylinder. In other words
$$
K_\zeta(\zeta_n^x) \leq K_\zeta(\zeta_{\tau(\zeta_n^x)}^x) + \log n.
$$
 From which follows
\begin{eqnarray*}
\underline K_\zeta(x,T)
&\leq &
\liminf_{n\to\infty}[K_\zeta(\zeta_{\tau(\zeta_n^x)}^x)+\log n]/n \\
&=&
\liminf_{n\to\infty}
\frac{K_\zeta(\zeta_{\tau(\zeta_n^x)}^x)}{\tau(\zeta_n^x)} \times \frac{\tau(\zeta_n^x)}n \\
&\leq&
\overline K_\zeta(x,T) \underline R_\zeta(x).
\end{eqnarray*}
White's  improvement of Brudno's theorem \cite{w,b} gives for $\mu$-a.e. $x\in X$
$$
\underline K_\zeta(x,T) = \overline K_\zeta(x,T) = h_\mu(T|\zeta) > 0,
$$
hence $\underline R_\zeta(x)\geq 1$.
\end{proof}
After we discovered this proof of Theorem~\ref{thm:tau/n} Saussol
gave an alternate proof using the Shannon McMillan Breiman theorem instead
of Brudno's result \cite{acs}.

\section{Dimension and Lyapunov exponent via return times}\label{sec:low}

\subsection{Preliminaries}

We can apply the results of the previous section
to a very general case of one-dimensional piecewise monotonic maps.
For a function $f: [0,1] \to \RR$ and $p > 0$ we define the $p$--variation 
of $f$ by 
$$\text{var}_p(f) :=\sup \left \{ \sum_{i=1}^{N-1} 
\left | f(x_{i+1}) - f(x_i) \right |^p  \right \},
$$
where the supremum is taken along all finite ordered sequences
of points $0 \le x_1 < x_2 < \cdots < x_N \le 1 $ and integers $N$.

Throughout this section $T\colon [0,1]\to [0,1]$ is a piecewise
monotonic transformation which preserves the ergodic invariant
measure $\mu$, and $\Z$ denotes the finite $\mu$-partition (i.e. 
partition modulo $\mu$)
of the interval into monotonic pieces. We say that a measurable
function $T':[0,1]\to\RR$ is a derivative of $T$ if
\[
\int_a^b\glu T'(x)\,dx = T(b)-T(a)
\]
for any interval $[a,b]$ contained in some element of $\Z$. 
We then denote the Lyapunov exponent of an invariant measure $\mu$ by 
$$\lambda_\mu = \int \log |T'|d\mu.$$
Let $\phi = \log |T'|$ and set $S_n\phi = \sum_{i=0}^{n-1}\phi\circ T^i$.
Given a $\mu$-partition $\Y$ we denote by $\Y_n = \vee_{i=0}^{n-1} T^{-i}\Y$
its refinement and we denote by $Y_n^x$ the unique element of $\Y_n$
containing $x$. Notice that such an element exists and is unique for $\mu$-a.e. $x$.

\subsection{Balls and cylinder sets}
In this section we slightly adapt the results by Hofbauer and Raith \cite{hr}
(see also \cite{h}) in order to get a good comparison between balls and cylinders.
We denote by $|J|$ the length of an interval $J \subset [0,1]$.

\begin{proposition}[\cite{hr}]\label{pro:1}
Let $T$ be a piecewise monotonic transformation with a derivative of
bounded $p$-variation for some $p>0$. Let $\mu$ be an ergodic 
$T$-invariant measure with Lyapunov exponent $\lambda_\mu>0$.
Then for any $\e>0$ we have
\begin{enumerate}
\item\label{*a}
there exists a finite or countable $\mu$-partition 
$\Y$  with finite entropy into intervals which refines $\Z$;
\item\label{*b}
the partition $\Y$ is a generator, in particular $h_\mu(T,\Y)=h_\mu(T)$;
\item\label{*c}
for any $n$ and $x$ we have
$$
\left|
S_n\phi(x) -\log \frac{1}{|Y_n^x|} \int_{Y_n^x} \exp(S_n\phi(y))dy
\right| \le n\e;
$$
\item\label{*d}
for $\mu$-almost every $x$ we have
\[
\lim_{n\to\infty}\frac {1}{n} \log d(T^nx,\partial T^n Y_{n}^x) = 0.
\]
\end{enumerate}
\end{proposition}
We can use this proposition and the technique described in \cite{h,hr}
to prove the following

\begin{lemma}\label{lem:1}
Let $\Y$ be the partition given by Proposition~\ref{pro:1}. 
For $\mu$-a.e. $x$ the set of
accumulation points of the sequence $ -\frac{1}{n} \log |Y_n^x|$
lies in the interval $[\lambda_\mu-\e,\lambda_\mu+\e]$.
\end{lemma}

\begin{proof}
We have $|T^n Y_n^x| = \int_{Y_n^x} \exp( S_n\phi(y) ) dy$, hence
Proposition~\ref{pro:1}.c gives 
\begin{equation}\label{*1}
\left| S_n\phi(x) - \log |T^n Y_n^x| + \log |Y_{n}^x| \right| \le n\e.
\end{equation}
Since $ d(T^n x,\partial T^n Y_{n}^x) \le  |T^n Y_{n}^x| \le 1$ 
Proposition~\ref{pro:1}.d also implies that
$\lim_{n\to\infty} \frac 1n \log |T^n Y_{n}^x| = 0$
for $\mu$-a.e. $x$. Furthermore the Birkhoff ergodic theorem gives that
$\lim_n \frac 1n S_n\phi = \lambda_\mu$ for $\mu$-a.e. $x$, thus
using \eqref{*1} we get the result. 
\end{proof}

\begin{lemma}\label{lem:2}
Let $\Y$ be the partition given by Proposition~\ref{pro:1}. 
For $\mu$-a.e. $x$ the set of
accumulation points of the sequence $-\frac{1}{n} \log d (x,\partial Y_{n}^x)$
lies in the interval $[\lambda_\mu-2\e,\lambda_\mu+2\e]$.
\end{lemma}

\begin{proof}
By the mean value theorem we have
\[
 d(x,\partial Y_{n}^x) \inf_{Y_n^x}| (T^n)'|
\le
d(T^nx,\partial T^n Y_{n}^x) 
\le d(x,\partial Y_{n}^x) \sup_{Y_n^x}| (T^n)'|.
\]
Since the logarithm is increasing this implies
\[
\inf_{Y_n^x} \log | (T^n)'| \le \log \frac{d(T^nx,\partial T^n Y_{n}^x)}{d(x,\partial Y_{n}^x)}
  \le \sup_{Y_n^x} \log | (T^n)'|
\]
Using \eqref{*1} this yields 
\begin{equation}\label{*2}
\begin{split}
\left|\log \frac{d (x,\partial Y_{n}^x)}{d(T^n x,\partial T^n Y_{n}^x) } + S_n\phi(x) \right|
&\le 2 \sup_{y\in \Y_n^x} |\phi(y)-\phi(x)|\\
&\le 2 n\e,
\end{split}
\end{equation}
by Proposition~\ref{pro:1}.c.
In addition, by the Birkhoff ergodic theorem we have 
$\lim_{n\to\infty} \frac 1nS_n\phi(x)=\lambda_\mu$ for $\mu$ a.e.~$x$. 
The conclusion follows then from Proposition~\ref{pro:1}.d
and inequality~\ref{*2}.
\end{proof}

\subsection{A lower bound for the Lyapunov exponent}

We are now ready to state and prove the following result.
\begin{theorem}\label{thm:lower}
Let $T$ be a piecewise monotonic transformation with a derivative
of bounded $p$-variation for some $p>0$. If $\mu$ is an ergodic $T$-invariant measure
with non-zero entropy, then
\begin{equation}\label{**1}
\lambda_\mu \ge \left (\liminf_{r\to 0} \frac{\tau(B_r(x))}{-\log r} \right )^{-1} 
\end{equation}
for $\mu$-almost every $x$.
\end{theorem}

\begin{remark}
1) Notice that each $C^{1+\e}$ piecewise monotonic map with finitely many pieces 
has a derivative of bounded $p$-variation, for $p\geq1/\e$, hence $C^{1+\e}$ multimodal
maps with non-zero entropy satisfies hypotheses of Theorem~\ref{thm:lower}.
\end{remark}

\begin{proof}
By Ruelle's inequality, the assumption that the entropy is positive
implies that the Lyapunov exponent is positive as well
\cite{h2}.
Let $\e>0$ and $\Y$ be the partition given by Proposition~\ref{pro:1}.

Let $x\in [0,1]$ be fixed.
We set $d_n = d(x,\partial Y_n^x)$ and $D_n = |Y_n^x|$.
Observe that since $\Y$ is generating we have $\lim_n D_n=0$,
hence $d_n$ converges monotonically to zero. 
Thus given $r>0$ we can define $n(r)$ to be the smallest 
integer $n$ such that $d_{n+1} < r \leq d_n$.
Note that we have $B_r(x) \subset Y_{n(r)}^x$, which implies
that $\tau( B_r(x) ) \ge \tau ( Y_{n(r)}^x )$.
Since $r\ge d_{n(r)+1}$ and $n(r)\to \infty$ as $r\to 0$ we get
\[
\begin{split}
\liminf_{r\to 0} \frac{\tau(B_r(x))}{-\log r} 
&\ge 
\liminf_{r\to 0}
\left(\frac{\tau(Y_{n(r)}^x)}{n(r)} \times 
\frac{1}{-\frac{1}{n(r)}\log d_{n(r)+1}}\right) \\
&\ge 
\left(\liminf_{n\to\infty} \frac{\tau(Y_{n}^x)}{n}\right)
\times
\left(\limsup_{n\to\infty} -\frac{1}{n} \log d_{n}\right)^{-1}.
\end{split}
\]
Since $h_\mu(T,\Y) = h_\mu(T) > 0$ by 
Proposition~\ref{pro:1}.b we can apply Theorem~\ref{thm:tau/n}, hence there exists a 
set $X_\e^1$ of full $\mu$-measure such that 
$\liminf_{n\to\infty} \frac{\tau(Y_{n}^x)}{n}\ge 1$ for any $x\in X_\e^1$.
By Lemma~\ref{lem:2} there exists a set of full
$\mu$-measure $X_\e^2$ such that 
$\limsup_{n\to\infty} -\frac{1}{n} \log d_{n} \le \lambda_\mu+2\e$
for any $x\in X_\e^2$. Thus for any $x\in X_{\e}^1 \cap  X_{\e}^2$ we get
\[
\liminf_{r\to 0} \frac{\tau(B_r(x))}{-\log r} \ge \frac{1}{\lambda_\mu+2\e}.
\]
We conclude that the inequality \eqref{**1} holds on the set of full measure
$\cap_{i\ge 1} (X_{1/i}^1 \cap  X_{1/i}^2)$. This proves the theorem. 
\end{proof}

In the case of Markov expanding maps we get a stronger result
\begin{corollary}
Under the hypotheses of Theorem~\ref{thm:lower}, if in addition
$T$ is piecewise expanding Markov then 
\[
\lim_{r\to 0} \frac{\tau(B_r(x))}{-\log r} = \frac{1}{\lambda_\mu}.
\]
\end{corollary}
\begin{proof}
Without loss of generality we assume that $T$ is topologically mixing.
We only sketch the proof.
If $\Z$ is a Markov partition for $T$ then it is easy to see
that for any $\e>0$ the partition $\Y = \Z_p = \vee_{i=0}^{p-1}\Z$
will have all the properties mentioned in Proposition \ref{pro:1},
provided $p$ is chosen sufficiently large.
Furthermore, $\Z_p$ is still a Markov partition, hence it has the 
specification property. Taking into account Remark~\ref{rem:1},
and proceeding as in the second part of the proof of Theorem~\ref{thm:dim}
yields to the conclusion.
\end{proof}

\subsection{Dimension via return time}

Give a map $T: X \to X$ on the metric space $(X,d)$ we 
define the first return of a point $x\in X$ into its $r$-ball $B_r(x)$ 
by
$$
\tau_r(x) := \min \left \{ k > 0 : T^k x \in B_r(x) \right \}
=
\min  \left \{ k>0 : d(T^k x,x) <r \right \}
$$
Given a measurable $\mu$-partition $\Y$ we denote by 
\[
R_n(x,\Y) = \inf \{k>0 : T^k x \in Y_n^x\}
\]
the repetition time of the first $n$ symbols of $x$. Ornstein and Weiss
have proven \cite{ow} that, whenever $\Y$ is a finite measurable $\mu$-partition
we have 
\begin{equation}\label{***1}
\lim_{n\to\infty} \frac{\log R_{n}(x,\Y)}{n} = h_\mu(T,\Y),
\end{equation}
for $\mu$-almost every $x$. See also \cite{ka,q} for the generalization to 
the case of a countable partition $\Y$.
This result will be essential to prove the following.
\begin{theorem}\label{thm:dim}
Let $T$ be a piecewise monotonic transformation with a derivative
of bounded $p$-variation for some $p>0$. If $\mu$ is an ergodic $T$-invariant measure
with non-zero entropy, then
\begin{equation}\label{**2}
d_\mu(x) = \lim_{r\to 0} \frac{\log\tau_r(x)}{-\log r}
\end{equation}
for $\mu$-almost every $x$.
\end{theorem}

\begin{proof}
By Ruelle's inequality, the assumption that the entropy is positive
implies that the Lyapunov exponent is positive as well
\cite{h2}.
Let $\e\in (0,\lambda_\mu)$ and $\Y$ be the partition given by Proposition~\ref{pro:1}.

We proceed as in the proof of Theorem~\ref{thm:lower},
and keep the same notations. 
Observe that for any $x$ and $r>0$ we have $B_r(x)\subset Y_{n(r)}^x$, from
which follows $\tau_r(x) \ge R_{n(r)}(x,\Y)$. Thus
\[
\begin{split}
\liminf_{r\to 0} \frac{\log\tau_r(x)}{-\log r}
&\ge
\liminf_{r\to 0} \left(
\frac {\log R_{n(r)}(x,\Y)}{n(r)} \times \frac {1}{-\frac{1}{n(r)} \log d_{n(r)+1}} 
\right)\\
&\ge
\left(\liminf_{n\to\infty} \frac{\log R_{n}(x,\Y)}{n}\right)
\times 
\left(\limsup_{n\to\infty} -\frac{1}{n} \log d_{n}\right)^{-1}.
\end{split}
\]
Since $h_\mu(T,\Y) = h_\mu > 0$ by Proposition~\ref{pro:1}.b 
we can apply the countable alphabet version of Ornstein and Weiss
return times theorem \cite{ow,ka,q}, hence there exists a set $X_\e^1$ 
of full $\mu$-measure such that \eqref{***1} holds for any $x\in X_\e^1$.
By Lemma~\ref{lem:2} there exists a set of full $\mu$-measure $X_\e^2$ such that 
$\limsup_{n\to\infty} -\frac{1}{n} \log d_{n} \le \lambda_\mu+2\e$
for any $x\in X_\e^2$. 
Thus for any $x\in X_{\e_i}^1 \cap  X_{\e_i}^2$ we get
\begin{equation}\label{12}
\liminf_{r\to 0} \frac{\log \tau_r(x)}{-\log r} \ge \frac{h_\mu}{\lambda_\mu+2\e}.
\end{equation}

Next we want to find an upper bound for the $\limsup$ of the same quantity.
If $m(r)$ denotes the smallest integer $m$ such that $D_{m} < r \le D_{m-1}$,
then we have $B_r(x) \supset Y_{m(r)}^x$. Thus
\[
\begin{split}
\limsup_{r\to 0} \frac{\log\tau_r(x)}{-\log r}
&\le
\limsup_{r\to 0} \left(
\frac {\log R_{m(r)}(x,\Y)}{m(r)} \times \frac {1}{-\frac{1}{m(r)} \log D_{m(r)-1}} 
\right)\\
&\le
\left(\limsup_{m\to\infty} \frac{\log R_{m}(x,\Y)}{m}\right)
\times 
\left(\liminf_{m\to\infty} -\frac{1}{m} \log D_{m}\right)^{-1}.
\end{split}
\]
By Lemma~\ref{lem:1} there exists a set of full $\mu$-measure $X_\e^3$
such that for any $x\in X_\e^3$ we have 
$\liminf_{m\to\infty} -\frac{1}{m} \log D_{m}  \ge \lambda_\mu-\e$.
This together with \eqref{***1} implies that for any $x\in X_\e^1 \cap X_\e^3$ we have
\begin{equation}\label{13}
\limsup_{r\to 0} \frac{\log \tau_r(x)}{-\log r} \le \frac{h_\mu}{\lambda_\mu-\e}.
\end{equation}
In addition, Hofbauer \cite{h} has shown in this setting that 
\[
d_\mu(x) = \frac{h_\mu}{\lambda_\mu},
\]
hence we conclude by \eqref{12} and \eqref{13} that the equality \eqref{**2} holds
on the set of full measure $\cap_{i\ge 1} (X_{1/i}^1 \cap  X_{1/i}^2 \cap X_{1/i}^3)$. 
This finishes the proof.
\end{proof}

\section{Appendix:  Kolmogorov complexity and Brudno's theorem}

The idea of Kolmogorov complexity is that a finite 0--1 word is 
only as complicated as the algorithm that produces it.\footnote{The 
generalization to arbitrary finite alphabets is straightforward.}  
To run an 
algorithm we need to fix a computer (with infinite storage capacity).
The Kolmogorov complexity $K_M(x)$ of a word $x$ with respect to a fixed 
computer $M$ is the length of the shortest algorithm which outputs $x$ given 
the length of $x$ as an input. 
Kolmogorov proved that there exist universal
computers $U$ such that 
\begin{equation}
K_U(x) \le K_M(x) + C
\label{optimal}
\end{equation}
where $C$ is a constant depending only on $U$ and $M.$  Here the word
universal is used to indicate that $U$ can simulate any other computer $M$.

More formally, a computer $M$ is a Turing machine while an algorithm 
which produces a finite 0--1 string $s$ is a 0--1 string $p$ such that $M(p)
= s$.  If there is no $p$ with $M(p) = s$ we say that the length $l(p)$ of the
algorithm is not defined while if there is more than one such $p$ we choose 
the first in the lexicographical order.

If $p$ is a finite word of length $n$ then we denote by 
$\hat{p}$ the string
$$p(0)p(0)p(1)p(1)\cdots p(n-1)p(n-1)01.$$  
If we input the concatenated word $\hat{p} q$ into a suitably programmed Turing
machine it will recognize two distinct inputs: $p$ and $q$.  Also, 
if $n \in \bf{N}$
let $[n]$ be the $n$ binary string in the lexicographical order given by
$$0, \ 1, \ 00, 01, \ 10, \ 11, \ 000, \ \dots$$ 
i.e.~$[3]=00$. Notice that
$l([n]) \le \log_2 n.$

There are countably many Turing machines, which may be computable enumerated
as $A_1,A_2,\dots$.  We say that a Turing machine is {\em universal} if,
for any $m$ and any finite word $p: U(\hat{[m]} p) = A_m(p)$. Thus a universal
Turing machine simulates any given machine on any given input.

For $x$ an infinite 0--1 sequence one defines the average complexity by 
looking at the first $n$--bits $x(n)$ and defining
$$\overline{K}(x) := \limsup_{n \to \infty} \frac{K_U(x(n))}{n}.$$ 
Note that by equation \eqref{optimal}
the average complexity does not depend which universal computer $U$
is chosen.
The function $\underline{K}(x)$ is defined in an analogous way with
the $\limsup$ replaced by a $\liminf$.

Brudno's theorem shows the linkage between
complexity and entropy.  Suppose that $\mu$ is an ergodic invariant measure
for the map $f$ and $\zeta$ is a finite measurable partition.  The partition
$\zeta$ gives rise to the symbolic space $\Sigma_{\zeta}$ and a map
$\phi:X \to \Sigma_{\zeta}$ which is a semiconjugacy $\sigma \phi = \phi T.$
Let $\zeta_n^x$ be the word consisting of the first $n$--symbols of the 
sequence $\phi(x)$ and define 
$\underline{\overline{K}}_{\zeta}(x,T) = \underline{\overline{K}}((\zeta_n^x)_n)$.
Brudno has shown:

\begin{theorem}
{\bf (Brudno \cite{b})} 
$\overline{K}_{\zeta}(x) = h_{\mu}(f|\zeta)$ for $\mu$--almost every point $x$
\end{theorem}

White has improved this theorem, he has shown:

\begin{theorem} {\bf (White \cite{w,w1})}
$\underline{K}_{\zeta}(x) = \overline{K}_{\zeta}(x) = h_{\mu}(f|\zeta)$ for $\mu$--a.e.~$x$.
\end{theorem}


\begin{thebibliography}{99}
\footnotesize{

\bibitem{a} Afraimovich, V.~{\it Pesin's dimension for Poincar\'e recurrence},
Chaos {\bf 7} (1997) 12--20.

\bibitem{acs} Afraimovich,~V., Chazottes,~J.-R., and Saussol,~B.~{\em
Pointwise dimensions for Poincar\'e recurrence associated with maps and special flows},
preprint 2000.

\bibitem{bs} Barreira,L.~and Saussol,~B.~{Hausdorff dimension of measures
via Poincar\'e recurrence}, Comm.~Math.~Phys.~{\bf 219} (2001) 443-464.
 
\bibitem{bs2}  Barreira,L.~and Saussol,~B.~{Product structure of Poincar\'e
recurrence}, to appear in ETDS.

\bibitem{b} Brudno,~A.~{\em Entropy and the complexity of the trajectories
of a dynamical system}, Russ.~Math.~Surv.~{\bf 2} (1983) 
127--51 (Engl.~Trans.).

\bibitem{chv} Cassaigne,~J., Hubert,~P.~and Vaienti,~S.~in preparation.

\bibitem{hv} Haydn,~N.~and Vaienti,~S.~{\it The limiting distributions and
error terms for return times of dynamical systems}, preprint (2001).

\bibitem{h} Hofbauer,~F.~{\em Local dimension for piecewise monotonic maps on the interval}, Erg.~Th.~Dyn.~Sys.~{\bf 15} (1995) 1119--1142.

\bibitem{h2}  Hofbauer,~F.~{\em An inequality for the 
Ljapunov exponent of an ergodic invariant measure for a 
piecewise monotonic map on the interval}, in: Lyapunov 
exponents, Proceedings, Oberwolfach, 1990 (Eds.: L. Arnold, 
H. Crauel, J.-P. Eckmann), Lecture Notes in Mathematics 1486, 
Springer, Berlin, 1991, pp. 227-231.

\bibitem{hr} Hofbauer,~F. and Raith,~P.~{\em The Hausdorff dimension of 
an ergodic invariant 
measure for a piecewise monotonic map of the interval}, 
Canad. Math. Bull.~{\bf 35} (1992) 84--98.

\bibitem{hsv} Hirata,~M., Saussol,~B.~and Vaienti,~S.~{\em Statistics of 
return times: a general framework and new applications}, 
Comm.~Math.~Phys.~{\bf 206} (1999) 33--55.

\bibitem{kh} Katok, T.~and Hasselblatt, B.~{\it Introduction to the modern
theory of dynamical systems}, Cambridge Univ.~Press (1995).

\bibitem{ka} Kontoyiannis,~I., Algoet,~P., Suhov,~Yu.~and Wyner,~A.~{\it 
Nonparametric entropy esitmation for stationary processes and random fileds,
with applications to English text}, IEEE Tras.~Inf.~Th.~{\bf 44} (1998)
1319--1327.
1988

\bibitem{lpv} Haydn,~N., Luevano,~J., Mantica, G.~and Vaienti, S.~{\it 
Multifracatal properties of return time statistics},
submitted to Phys Rev Letters, 2001

\bibitem{ow} Ornstein,~D.~and Weiss,~B.~{\it Entropy and data compression},
IEEE Trans. Inf.~Th.~{\bf 39} (1993) 78--83. 

\bibitem{psv} Penn\'e,~V., Saussol,~B.~and Vaienti,~S.~{\it Dimensions for 
recurrence times: topological and dynamical properties}, 
Disc.~Cont.~Dyn.~Sys.~{\bf 4} (1998) 783--798.

\bibitem{q} Quas,~A.~{\it An entropy esitmator for a class of infinite
alphabet processes}, Theor.~Veroyatnost.~i Primenen.~{\bf 43} (1998)
61---621.

\bibitem{w} White,~H. {\em Algorithmic complexity of points in a dynamical
system}, Erg. Th.~Dyn.~Sys.~{\bf 13} (1993) 807--30.

\bibitem{w1} White,~H. {\em On the algorithmic complexity of trajectories of 
points in dynamical systems}, Ph.D.~dissertation Univ.~of North Carolina 
at Chapel Hill 1991.}

\end{thebibliography}
\end{document}